\documentclass[12pt]{amsart}
\usepackage{amssymb,amsmath}
\usepackage{enumerate}
\usepackage{amsthm}
\usepackage[dvips]{color}
\usepackage[dvips]{graphicx}
\usepackage[dvips]{graphics}
\usepackage{keyval}
\usepackage{varioref}
\usepackage{amsfonts}
\usepackage[mathscr,mathcal]{eucal}
\usepackage{longtable}

\theoremstyle{plain}
\newtheorem{thm}{Theorem}[section]
\newtheorem{cor}[thm]{Corollary}
\newtheorem{lem}[thm]{Lemma}
\newtheorem{prop}[thm]{Proposition}
\newcommand{\sig}{\sigma}

\newcommand{\ld}{\vartriangleleft}

\newcommand{\f}[1]{\mathfrak{#1}}

\newcommand{\N}{\mathbb{N}}
\newcommand{\Z}{\mathbb{Z}}
\newcommand{\Bo}{$\square$}
\DeclareMathOperator{\Aut}{Aut} \DeclareMathOperator{\Inn}{Inn}
\DeclareMathOperator{\Core}{Core}
\DeclareMathOperator{\Stab}{Stab} 

\def\Rem{{\bf Remark.}\;}

\begin{document}
\title{The STRUCTURE OF $G/\Phi(G)$ IN TERMS OF $\sigma(G)$}
\author{Alireza Jamali\; \&\; Hamid Mousavi}
\address{Institute of Mathematics, University for Teacher Education,
 599 Taleghani Avenue, Tehran 15614, Iran}
\email{jamali@saba.tmu.ac.ir}
 \address{ Institute for Advanced Studies in Basic Sciences, P.O.Bax 45195-159, Gavazang, Zanjan, Iran}
 \email{hmousavi@iasbs.ac.ir}
\date{}
\begin{abstract}
Let $G$ be a finite group. We let $\f{m}(G)$ and $\sig(G)$ denote
the number of maximal subgroups of $G$ and the least positive
integer $n$ such that $G$ is written as the union of $n$ proper
subgroups, respectively. In this paper we determine the structure
of $G/\Phi(G)$ when $G$ is a finite soluble group with
$\f{m}(G)\leq 2\sig(G)$.
\end{abstract}
\maketitle
\section{Introduction}
Let $G$ be a finite group. We let $\f{m}(G)$ and $\sig(G)$ denote
the number of maximal subgroups of $G$ and the least positive
integer $n$ such that $G$ is written as the union of $n$ proper
subgroups, respectively.

In [\ref{JaMo}], the authors produce a number of results which
enables us to determine the structure of $G/\Phi(G)$ when $G$ is a
finite soluble group with $\f{m}(G)\leq 2\sig(G)$. The
determination of $G/\Phi(G)$when $G$ is soluble is of particular
interest, for the Hall subgroups of $G/\Phi(G)$ gives us some
information about the Hall subgroups of $G$. Pazderski [\ref{Paz}]
proves that the groups having at most $20$ maximal subgroups are
soluble. The author in [\ref{Wang}] gives, among the others, a
result which simplifies the proof of  Pazderski's result as well
as classifying those finite groups that have exactly eight maximal
subgroups. In [\ref{Thes}], the author has determined the
structure of $G/\Phi(G)$, where $G$ is a finite non-nilpotent
group with $\f{m}(G)\leq 20$. For an application, the author
[\ref{Thes}] gives, for example, a complete list of finite
non-nilpotent groups having exactly 9 maximal subgroups.

The authors [\ref{JaMo}] proves the following results:

\begin{thm}\label{1.1}
Let $G$ be a non-cyclic finite soluble group with $\sig(G)=n$.
Then there exists $n$ maximal subgroups $M_1,$...,$M_n$ of indices
$i_1,$...,$i_n$, respectively, such that
$G=\bigcup^n_{\ell=1}{M_{\ell}}$,
$i_1\leq\cdots\leq i_n$ and one of the following statements holds:\\
1) $i_1< n-1$ and $M_2,$...,$M_n$ are conjugate;\\
2) $i_1=n-1$ and $M_i\lhd G$ for all $i=1,2,\dots ,n$.
\end{thm}

\vspace*{0.2cm}An $n$-tuple $(M_1,..., M_n)$ of maximal subgroups
of $G$ is called a $\sig$-cover of $G$ whenever
$G=\bigcup^n_{\ell=1}{M_{\ell}}$. A $\sig$-cover of $G$ satisfying
the condition (1) of above result is said to be a {\it conjugate
$\sig$-cover} of $G$ while that satisfying the condition (2) is
called a {\it normal $\sig$-cover} of $G$.

\begin{prop}\label{1.4}
Let $G$ be a finite soluble group. If $G$ has a non-normal maximal
subgroup of index $\sig(G)-1$, then $G$ has a conjugate
$\sig$-cover.
\end{prop}

\begin{lem} \label{1.2} Let $G$ be a finite soluble group, and let $p$
be a prime divisor of $|G|$. If $G$ has at least two maximal
subgroups of index $p^n$, then $\sig(G)\leq 1+p^n$.
\end{lem}

\begin{cor}\label{1.3}
Let $G$ be a soluble group with $\sig(G)=n$. If $m<n-1$, then $G$
has at most one maximal subgroup of index $m$.
\end{cor}

In this paper, by using the above results, for a finite
non-nilpotent soluble group $G$, we give the structure of
$G/\Phi(G)$ when $G$ has at most $2\sig(G)$ maximal subgroups.

\section{The structure of $G/\Phi(G)$}
In this section we aim to obtain the structure of $G/\Phi(G)$,
where $G$ is a non-nilpotent soluble group. We shall show that the
structure of $G/\Phi(G)$ depends completely on the value of
$\sig(G)$.

\vspace*{0.3cm} For simplicity we shall use the abbreviations
$\f{m}$ and $\sig$ for $\f{m}(G)$ and $\sig(G)$, respectively.
Clearly $\sig\leq\f{m}$. Following Cohn [\ref{Cohn}], a group $G$
is said to be a {\it primitive $n$-sum group} if $G$ has no normal
subgroup $N$ for which $\sig(G/N)=\sig(G)$.

\begin{prop}\label{3.1}
Let $G$ be any finite group with $\f{m} >\sig$, and let
$(M_1,\cdots ,M_{\sig})$ be a $\sig$-cover of $G$. Suppose that
$M_i\ld G$ for $\sig<i\leq \f{m}$, and that
$H=\bigcap_{i=1}^{\sig}M_i$. Then $H\ld G$ and $G/\Phi(G)\cong
H/\Phi(G)\times G/H$. Moreover,
\begin{itemize}
  \item [{\rm(i)}] $H/\Phi(G)$ is a direct product of elementary
  abelian groups;
  \item  [{\rm(ii)}] if $G$ has a unique $\sig$-cover, then $G/H$
  is a primitive $\sig$-sum group.
\end{itemize}
\end{prop}

\noindent {\bf Proof.} The normality of $H$ is obvious. We set
$H_0=H$, $l_0=\sig+1$ and for $j\geq 1$, we define
\begin{align*}
l_j&=\min\{i\; | \; H_{j-1}\not\leq M_i, l_{j-1}\leq i\leq\f{m}\},
\\ H_j&=H_{j-1}\cap M_{l_j}.
\end{align*}
Hence $\Phi(G)=H_n$ for some integer $n$. Now on taking
$K_j=\bigcap_{i=1}^{j}M_{l_i}$\break $(1\leq j\leq n)$, we find
that $G=H_{j-1}M_{l_j}$ $(1\leq j\leq n)$, and
$G=K_{j}M_{l_{j+1}}$ \mbox{$(1\leq j\leq n-1)$.} Therefore,
$|H_n||G|^n=|H_0|\prod_{j=1}^n |M_{l_j}|$, and
$|K_n||G|^{n-1}=$\break $|K_1|\prod_{j=2}^n |M_{l_j}|$. Since
$K_1=M_{l_1}$ and $H_n=H_0\cap K_n$, we get $|G|=|H_0K_n|$. So
$G=HK_n$ and $G/H\cong K_n/\Phi(G)$. This proves the main part of
the proposition, as $G/\Phi(G)\cong H/\Phi(G)\times K_n/\Phi(G)$.
Now since $H/\Phi(G)\cong G/K_n$ and $G/K_n$ has no non-normal
maximal subgroup, because $K_n\nleqslant M_j$ for $j\leq \sigma$,
the group $G/K_n$ is nilpotent; but $\Phi(G/K_n)$ is trivial and
hence (i) is proved.

\vspace*{0.2cm} For the second part of the proposition, assume
that $N/H$ is a non-trivial normal subgroup of $G/H$ and that
$\sig(G/N)=\sig(G/H)$. Let $(M^*_1/N,\cdots ,M^*_{\sig}/N)$ be a
$\sig$-cover of $G/N$. Clearly $(M^*_1,\cdots ,M^*_{\sig})$ is a
$\sig$-cover of $G$. By the uniqueness assumption, we then have
$\bigcap_{i=1}^{\sig}{M^*_i}=H$. Therefore, $H=N$, a
contradiction. \Bo

\begin{cor}\label{3.1.5}
If a finite soluble group $G$ has a non-normal abelian maximal
subgroup of index $\sig(G)-1$, then $G/\Phi(G)\cong \Inn(G)\times
Z(G)/\Phi(G)$ and $\Inn(G)$ is primitive $\sig$-sum group.
\end{cor}

\noindent {\bf Proof.}  Using the notation of
Proposition~\ref{3.1}, we find that $H=Z(G)$ by [3, Corollary
4.7]. Also $G$ has a unique conjugate $\sig$-cover, so that
$\Inn(G)$ is primitive $\sig$-sum group. \Bo

\vspace*{0.3cm} \noindent \Rem We note that if $\sig =\f{m}$ for a
finite group $G$, then $G/\Phi(G)$ is a primitive $\sig$-sum
group. Hence, if $\f{m}=3$ then $G/\Phi(G)\cong \Z_2\times\Z_2$,
by [\ref{Cohn}, Theorem~3]. Therefore $G$ is a $2$-group.

\vspace*{0.2cm}In what follows $E_n$ {\it denotes an elementary
abelian group of order $n$}.

\begin{lem}\label{3.2}
Let $G$ be a finite non-nilpotent soluble group. If $\sig=\f{m}$,
then $G/\Phi(G)\cong E_{\sig-1}\rtimes \Z_{p^{\alpha}}$,\; $\alpha
>0$, $\sig\neq 3$ and $\Z_{p^{\alpha}}$ acts faithfully on
$E_{\sig-1}$.
\end{lem}

\noindent {\bf Proof.} According to Theorem~\ref{1.1}, for each
$\sig$-cover $(M_1,\cdots ,M_{\sig})$ of $G$, we have $i_1<\sig
-1$ and $M_2,...,M_{\sig}$ are conjugate in $G$. Since $G/\Phi(G)$
is a primitive $\sig$-sum group, $G/\Phi(G)$ is primitive soluble
group. Therefore, by [\ref{DoHa}, Theorem~15.6],
$G/\Phi(G)=L/\Phi(G)\rtimes M/\Phi(G)$ where $L/\Phi(G)$ is an
elementary abelian group of order $\sig-1$, and $M/\Phi(G)$ is a
cyclic $p$-group. Clearly the action of $M/\Phi(G)$ on $L/\Phi(G)$
is faithful. \Bo

\begin{lem}\label{3.3}
Let $G$ be a finite non-nilpotent soluble group. Suppose that
$\sig<\f{m}$ and $(M_1,\cdots ,M_{\sig})$ is a $\sig$-cover of $G$
such that $M_i\ld G$ for each $i$, $(\sig<i\leq \f{m})$. Then
$G/\Phi(G)\cong(E_{\sig-1}\rtimes\Z_t)\times\prod_{i=1}^{\ell}{\Z_{p_i}}$,
where $t\big| |\Aut(E_{\sig-1})|$, $p_1,...,p_{\ell}$ are primes
(not necessarily distinct), $\sig >3$ and
$\ell=\f{m}(G)-\f{m}(\Z_t)-\sig+1$.
\end{lem}

\noindent {\bf Proof.} As $G$ is non-nilpotent, the $\sig$-cover
$(M_1,\cdots ,M_{\sig})$ of $G$ is unique. Now the lemma follows
from Proposition~\ref{3.1}, by observing that $G/H\cong
E_{\sig-1}\rtimes\Z_t$ and $H/\Phi(G)\cong
\prod_{i=1}^{\ell}{\Z_{p_i}}$ where $H=\bigcap_{i=1}^{\sig}M_i$.
\Bo

\begin{cor}\label{3.4}
Let $G$ be a finite non-nilpotent soluble group. If
$\f{m}<2\sig-1$ then $G/\Phi(G)\cong
(E_{\sig-1}\rtimes\Z_t)\times\Z_{p_1\cdots p_{\ell}}$, where
$t\big| |\Aut(E_{\sig-1})|$,\; $p_1,...,p_{\ell}$ are distinct
primes coprime to $t$, $\sig\neq 3$ and
$\ell=\f{m}(G)-\f{m}(\Z_t)-\sig+1$.
\end{cor}

\noindent {\bf Proof.} Suppose $(M_1,\cdots ,M_{\sig})$ is any
$\sig$-cover of $G$. By Corollary~\ref{1.3} it follows that each
$M_i$ with $\sig<i\leq m$ is unique. Therefore, $M_i\ld G$ for
$\sig<i\leq m$ and the result is proved by using Lemma~\ref{3.3}.
Now since the maximal subgroups of $H/\Phi(G)$, where
$H=\bigcap_{i=1}^{\sig}M_i$, are unique, we see that
$p_1,...,p_{\ell}$ are distinct. \Bo

\begin{thm}\label{3.4'}
Let $G$ be a soluble group having $m$ non-conjugate maximal
subgroups $M_1,...,M_m$. Suppose that $M_1,...,M_m$ are non-normal
in $G$. Let $C_i=\Core_G(M_i)$ for $1\leq i\leq m$ and
$H_i=\bigcap_{\underset {\ell\neq i}{\ell=1}}^{n}{C_{\ell}}$,
where $n\leq m$ and $1\leq i\leq n$.
\begin{itemize}
  \item [{\rm(i)}] If $H_i\nsubseteq C_i$ for all $i$, then
$G/\bigcap_{i=1}^{n}{C_i}\cong(\prod_{i=1}^{n}{L_i/C_i})
  \rtimes(\bigcap_{i=1}^{n}{M_i}/\bigcap_{i=1}^{n}{C_i})$,
  where $L_i$ is the normal minimal subgroup of $G$ containing $C_i$.
  \item [{\rm(ii)}] If $G=H_{\ell}C_{\ell}$ for some $\ell$, then
  $\bigcap_{i=1}^{n}{M_i}/\bigcap_{i=1}^{n}{C_i}\cong M_{\ell}/C_{\ell}
  \times \bigcap_{\substack{i=1\\i\neq \ell}}^{n}{M_i}/H_{\ell}$.
  \item [{\rm(iii)}] If $H_{\ell}\cap M_{\ell}\leq C_{\ell}$ for some
  $\ell$, then
  $\bigcap_{i=1}^{n}{M_i}/\bigcap_{i=1}^{n}{C_i}\cong \bigcap_{\underset{i\neq\ell
  }{i=1}}^{n}{M_i}/\bigcap_{\underset{i\neq\ell}{i=1}}^{n}{C_i}$.
\end{itemize}
\end{thm}

\noindent {\bf Proof.} For simplicity we set
$L=\bigcap_{i=1}^{n}{L_i}$, $C=\bigcap_{i=1}^{n}{C_i}$,
$M=\bigcap_{i=1}^{n}{M_i}$ and
$K_{\ell}=\bigcap_{\substack{i=1\\i\neq \ell}}^{n}{M_i}$.

(i) Since $L_i$ is unique, we have
\begin{equation*}
(L_i\cap H_i)M_i=(L_i\cap H_i)C_iM_i=(L_i\cap
H_iC_i)M_i=L_iM_i=G\quad (1\leq i\leq n).
\end{equation*}
\begin{equation*}
(L_i\cap H_i)/\bigcap_{\ell=1}^{n}{C_{\ell}}\cong (L_i\cap
H_i)C_i/C_i=(L_i\cap H_iC_i)/C_i=L_i/C_i\quad (1\leq i\leq n).
\end{equation*}
Now using $\bigcap_{i=1}^{n}{L_i}=(L_1\cap H_1)\cdots(L_n\cap
H_n)$, we obtain
\begin{align*}
L/C=(L_1\cap H_1)\cdots(L_n\cap H_n)/C
\cong\prod_{i=1}^{n}{\bigl((L_i\cap H_i)/C\bigl)}
\cong\prod_{i=1}^{n}{L_i/C_i},
\end{align*}
and $LM=(L_1\cap H_1)\cdots(L_n\cap
H_n)M=\bigcap_{i=1}^{n}{(L_i\cap H_i)M_i}=G.$ Therefore
\begin{align*}
  G\Bigl/\bigcap_{i=1}^{n}{C_i}=LM/C\cong
 L/C\rtimes M/C\cong (\prod_{i=1}^{n}{L_i/C_i})\rtimes
 \bigl(\bigcap_{i=1}^{n}{M_i}\Bigl/\bigcap_{i=1}^{n}{C_i}\bigl).
\end{align*}

\vspace*{0.2cm}(ii)We have the following isomorphism:
$$(H_{\ell}\cap M_{\ell})/C\cong (H_{\ell}\cap
M_{\ell})C_{\ell}/C_{\ell}=(H_{\ell}C_{\ell}\cap
M_{\ell})/C_{\ell}=M_{\ell}/C_{\ell}.$$ Now since $M=(M_{\ell}\cap
H_{\ell})(K_{\ell}\cap C_{\ell})$, we get
\begin{equation*}
 M/C=(M_{\ell}\cap H_{\ell})(K_{\ell}\cap C_{\ell})/C\cong
 (M_{\ell}\cap H_{\ell})/C\times (K_{\ell}\cap C_{\ell})/C\cong
 M_{\ell}/C_{\ell}\times K_{\ell}/H_{\ell}.
\end{equation*}

\vspace*{0.2cm}(iii) By our assumption, we find that $H_{\ell}\cap
M_{\ell}=\bigcap_{i=1}^{n}{C_i}$ and $H_{\ell}M_{\ell}=G$, whence
$$M/C\cong (K_{\ell}\cap M_{\ell})/(H_{\ell}\cap M_{\ell})\cong
(K_{\ell}\cap H_{\ell}M_{\ell})/H_{\ell}=K_{\ell}/H_{\ell}.\,
\square$$

\begin{cor}\label{cor}
With the above notation and assumption, if $|M_{\ell}/C_{\ell}|$
is a prime for some $\ell$ and if $ H_{\ell}C_{\ell}\neq G$, then
$\bigcap_{i=1}^{n}{M_i}/\bigcap_{i=1}^{n}{C_i}\cong(\bigcap_{\underset{i\neq
\ell}{ i=1}}^{n}{M_i})/H_{\ell}$.
\end{cor}

\noindent {\bf Proof.} We have $(H_{\ell}\cap
M_{\ell})/(H_{\ell}\cap C_{\ell})\cong (H_{\ell}C_{\ell}\cap
M_{\ell})/C_{\ell}\ld M_{\ell}/C_{\ell}$. But since
$M_{\ell}\not\ld G$, we conclude that $H_{\ell}\cap
M_{\ell}=H_{\ell}\cap C_{\ell}$. The result is now immediate by
using Theorem~\ref{3.4'} (part(iii)). \Bo

\vspace*{0.2cm}\noindent \Rem By (i) and (ii) of
Theorem~\ref{3.4'}, we have
\begin{equation}\label{1}
 G\Bigl/\bigcap_{i=1}^{n}{C_i}\cong
 \bigl(\prod_{i=1}^{n}{L_i/C_i}\bigl)\rtimes\bigl(\prod_{i=1}^{n}{M_i/C_i}\bigl).
\end{equation}
On the other hand, since $G=H_{\ell}C_{\ell}$, it is readily seen
that
\begin{equation}\label{2}
 G\Bigl/\bigcap_{i=1}^{n}{C_i}\cong G/C_{\ell}\times G/H_{\ell}.
\end{equation}
In what follows, we shall show that how the direct decomposition
of (\ref{2}) is obtainable from the semidirect decomposition of
(\ref{1}). To see this it is enough to show that $[L_i/C_i,
M_j/C_j]=1$ for all $i,$ $j$ with $i\neq j$. Since $M_j\cap H_j\ld
M_j$ and $L_i\cap H_i\ld M_j$ (for $i\neq j$), we have $$[L_i\cap
H_i, M_j\cap H_j]\leq (L_i\cap H_i)\cap (M_j\cap
H_j)=\bigcap_{\ell=1}^{n}{C_{\ell}}\quad (i\neq j),$$ and so
$\bigl[(L_i\cap H_i)/\bigcap_{\ell=1}^{n}{C_{\ell}}\;,(M_j\cap
H_j)/\bigcap_{\ell=1}^{n}{C_{\ell}}\bigr]=1$, for all $i, j$ with
$i\neq j$. Our assertion follows at once by considering the
following isomorphisms: $$L_i/C_i\cong (L_i\cap
H_i)\Bigl/\bigcap_{\ell=1}^{n}{C_{\ell}}\qquad ,\qquad
M_i/C_i\cong (M_i\cap H_i)\Bigl/\bigcap_{\ell=1}^{n}{C_{\ell}}.$$

\begin{thm}\label{3.5}
Let $G$ be a soluble group having two non-conjugate maximal
subgroups $M_1$ and $M_2$ such that $|G:C_1C_2|=\ell>1$ and
$M_i/C_i$ is cyclic of order $r_i$, $(i=1,2)$. Then $(M_1\cap
M_2)/(C_1\cap C_2)\cong \Z_t\times \Z_n$ where $t=r_1r_2/\ell n$
and $n=(r_1/\ell ,r_2/\ell)$.
\end{thm}

\noindent {\bf Proof.} Since $M_1$ and $M_2$ are non-conjugate
maximal subgroups of $G$, the subgroup $C_1C_2$ contains properly
both $C_1$ and $C_2$. Hence $L_i\leq C_1C_2\;(i=1,2)$, because
$L_i$ is the unique minimal normal subgroup of $G$ containing
$C_i$. Hence by $G/C_i\cong L_i/C_i\rtimes M_i/C_i \; (i=1,2)$;
we conclude that $\ell$ divides both $r_1$ and $r_2$. We now
distinguish two cases.

\vspace*{0.2cm} \noindent {\bf Case 1.} We assume that either
$C_1\cap M_2=C_1\cap C_2$ or $C_2\cap M_1=C_1\cap C_2$. If
$C_1\cap M_2=C_1\cap C_2$, then $(M_1\cap M_2)/(C_1\cap
C_2)=(M_1\cap M_2)/(C_1\cap M_2)=M_1/C_1$. Also $C_1C_2\cap
M_2=C_2(C_1\cap M_2)=C_2$. Hence $$r_2=|M_2:C_2|=|M_2:C_1C_2\cap
M_2|=|C_1C_2M_2:C_1C_2|=|G:C_1C_2|=\ell,$$ it follows that $t=r_1$
and $n=1$, and (ii) is proved at once. A similar argument can be
used when $C_2\cap M_1=C_1\cap C_2$.

\vspace*{0.2cm} \noindent {\bf Case 2.} We set $H_1=(C_1\cap
M_2)/(C_1\cap C_2)$ and $H_2=(C_2\cap M_1)/(C_1\cap C_2)$. Then
$H_1$ and $H_2$ are non-trivial. We let $K=(M_1\cap M_2)/(C_1\cap
C_2)$. Now $K$ is abelian, by using the fact that $(M_1\cap
M_2)/(C_1\cap C_2)\hookrightarrow M_1/C_1\times M_2/C_2$. We then
have $|K|=r_1r_2/\ell$, $|H_i|=r_j/\ell$, where $i\neq j$ and $i$,
$j\in\{1,2\}$. Clearly $K/H_i\cong \Z_{r_i}\; (i=1,2)$. So we may
assume that $K$ contains two elements $x$ and $y$ of order $r_1$
and $r_2$, respectively. On setting $H=\langle xy\rangle$, we get
$|H|=r_1r_2/\ell n$, where $n=(r_1/\ell ,r_2/\ell)$. Obviously, if
$n=1$, then $H=K$ and the proof is completed. We now suppose that
$n\neq 1$. Since $H_1\cap H_2=1$ and $H_1H_2$ is non-cyclic, we
conclude that $H$ cannot contain both $H_1$ and $H_2$. We assume
that $H_1\nleq H$, say. Then $K=HH_1$ and $|H\cap H_1|=r_2/n\ell$.
It is now straightforward to see that $K\cong \Z_t\times\Z_n$,
where $t=r_1r_2/n\ell$. \Bo

\begin{thm}\label{3.6}
Let $G$ be a finite non-nilpotent soluble group with $\f{m}=$
$2\sig-1$. Then $\sig\neq 3$ and
\begin{itemize}
\item[\rm{(i)}] $G/\Phi(G)\cong(E_{\sig-1}\rtimes
\Z_t)\times\Z_{p_1\cdots p_{\ell}}$, where $t\big |
|\Aut(E_{\sig-1})|$, $p_1,...,p_{\ell}$ are distinct primes
coprime to $t$ and $\ell=\sig-\f{m}(\Z_t)$. \item[\rm{(ii)}]
$G/\Phi(G)\cong(E_{\sig-1}\times E_{\sig-1})
\rtimes\Z_{p^{\alpha}}$, where $p$ is a prime and $\alpha\in\N$
such that $p^{\alpha}\big | |\Aut(E_{\sig-1})|$.
\end{itemize}
\end{thm}

\noindent {\bf Proof.} We first show that $G$ has a conjugate
$\sig$-cover. Suppose that $G$ has no conjugate $\sig$-cover and
assume that $M_1,...,M_{2\sig-1}$ are all maximal subgroups of $G$
such that $(M_1,\cdots ,M_{\sig})$ is a normal $\sig$-cover of
$G$. It follows from Corollary\ref{1.3} that,
$M_{\sig+1},...,M_{2\sig-1}$ are mutually conjugate in $G$. Since
$\sig-1$ is a prime, $G/C\cong \Z_{\sig-1}\rtimes(M_{\sig+1}/C)$,
where $C=\Core_G(M_{\sig+1})$. Hence, $G$ must have a maximal
normal subgroup of index less than $\sig-1$, which is impossible.
Therefore, $G$ has a $\sig$-cover $(M_1,\cdots ,M_{\sig})$
satisfying the condition (i) of the theorem~\ref{1.1}.

Now if $M_i\ld G$ for each $i>\sig$, then the first part of the
theorem occurs, and we have
$\ell=\f{m}(G)-\f{m}(\Z_t)+1=\sig-\f{m}(\Z_t)$.

Next suppose that $M_i\not\ld G$ for some $i>\sig$. Then
$|G:M_i|\geq\sig-1$ and hence $M_{\sig+1},...,M_{2\sig-1}$ form a
single conjugacy class in $G$. We set $C_1=\Core_G(M_2)$ and
$C_2=\Core_G(M_{\sig+1})$. Since $G$ has only one maximal subgroup
which is normal, $G\neq C_1C_2$. Also $M_2/C_1$ and
$M_{\sig+1}/C_2$ are cyclic, by Proposition 1.4(i). Now since
$M_1$ is the only maximal subgroup of $G$ which is normal, there
is a prime number $p$ such that $M_2/C_1$ and $M_{\sig+2}/C_2$
have $p$ power orders. By Theorem~\ref{3.5} the proof is
completed, where $p^{\alpha}=\max \{
|M_2/C_1|,|M_{\sig+1}/C_2|\}$. \Bo

\begin{lem}\label{3.7}
Let $q=2^n-1$ be a prime. Then there is a unique group of the
structure $\Z_2^n\rtimes\Z_q$.
\end{lem}

\noindent {\bf Proof.} To see this, it is enough to see that the
image of the corresponding action is a Sylow $q$-subgroup of
$\Aut(\Z_2^n)$. \Bo

\begin{prop}\label{3.8}
Let $G=LM$, where $L\ld G$, $L\cap M=1$ and $M$ be a maximal
subgroup of $G$. Suppose that the action of $M$ on $L\setminus
\{1\}$ is faithful and fixed-point-free. If $|L|-1$ is a prime and
$M$ is a Hamiltonian group, then
\begin{itemize}
\item[\rm(i)] $M$ acts transitively on $L\setminus\{1\}$;
\item[\rm(ii)] either $G\cong S_3$ or $G\cong \Z_2^n\rtimes\Z_q$,
where $q=|L|-1$.
\end{itemize}
\end{prop}

\noindent {\bf Proof.} Let $X$ be an arbitrary $M$-orbit of
$L\setminus\{1\}$ and let $H=\langle X\rangle$. Since $M$ is
maximal in $G$, $H=L$. It follows that $\Stab_M(X)$ is contained
in the kernel of the action, and hence $|X|=|M|$, because
$\Stab_M(x)=\Stab_M(X)$ for all $x\in X$. This shows that any
$M$-orbit of $L\setminus\{1\}$ has length equal to $|M|$.
Therefore, $|M|\big |(|L|-1)$ and we have $|M|=|L|-1$; that is,
$M$ acts transitively on $L\setminus\{1\}$, which proves (i). Now
by (i), we conclude that the elements of $L$ have the same prime
order $p$. So $|L|=p^n$ for some $n\in\N$. Since $|L|-1$ is a
prime, either $|L|=3$ or $p=2$. This completes the proof of (ii).
\Bo

\begin{lem}\label{3.9}
Suppose that $p,q$ are primes and $n,m$ are positive integers such
that $p^n=q^m+1$. Then one of the following statements holds
\begin{itemize}
\item[\rm(i)] $m=1$, and $q$ is a Mersenne prime;
\item[\rm(ii)] $n=1$, and $p$ is a Fermat prime;
\item[\rm(iii)] $q=n=2$ and $p=m=3$.
\end{itemize}
\end{lem}

\noindent {\bf Proof.} The proof is elementary. \Bo

\begin{thm}\label{3.10}
Let $G$ be a finite non-nilpotent soluble group with
$\f{m}(G)=2\sig$. Then $G/\Phi(G)$ is isomorphic to one of the
following groups:
\begin{itemize}
\item[\rm(i)] $\Z_2\times S_3$;
\item[\rm(ii)] $E_{\sig}\rtimes\Z_q\times\Z_q$, where $\sig
-1=q$ is a Mersenne prime;
\item[\rm(iii)] $E_{\sig -1}\rtimes\Z_t\times\Z_{p_1\cdots
p_{\ell}}$, where $p_i$'s are distinct primes being coprimes to
$t$, $\sig\neq 3$ and $\ell=\sig-\f{m}(\Z_t)+1$;
\item[\rm(iv)] $\Z_{\sig-1}\rtimes \Z_t\times\Z_{\sig-1}\times\Z_{\sig-1}$,
 where $\sig-1$ is a odd prime and $\f{m}(\Z_t)=1$;
\item[\rm(v)] $E_{\sig-1}\rtimes\Z_{t_1}\times
E_{\sig-1}\rtimes\Z_{t_2}$, where $t_1$ and $t_2$ are coprime
positive integers with $\f{m}(\Z_{t_1})=\f{m}(\Z_{t_2})=1$ and
$\sig-1$ is not a Fermat prime number; \item[\rm(vi)]
$(E_{\sig-1}\times E_{\sig-1})\rtimes\Z_t$, where $\f{m}(\Z_t)=2$
and $\sig\neq 3$; \item[\rm(vii)] $(E_{\sig-1}\times
E_{\sig-1})\rtimes\Z_t\times\Z_p$, where $\f{m}(\Z_t)=1$,
$(p,t)=1$, $t>2$ and $\sig\neq 3$; \item[\rm(viii)]
$E_{\sig}\rtimes(\Z_q \rtimes\Z_t)$, where $\sig-1=q$ is a
Mersenne prime and $\f{m}(\Z_t)=1$ with $t\big| \sig-2$.
\end{itemize}
\end{thm}

\noindent {\bf Proof.} Let $(M_1,\cdots ,M_{\sig})$ be a
$\sig$-cover of $G$, and let $M_{\sig+1},...,M_{2\sig}$ be the
remaining maximal subgroups of $G$. We set $C_i=\Core_G(M_i)$,
$1\leq i\leq 2\sig$. The proof is divided into two steps:

\vspace*{0.2cm} \noindent {\bf\it Step 1.}\\ We first assume that
$G$ has no conjugate cover. By Theorem~\ref{1.1}, we have $M_i\ld
G,\; 1\leq i\leq\sig$. Therefore $\sig-1$ is a prime. Since $G$ is
non-nilpotent, $M_{\sig+1}$ is non-normal with index $\sig$ in
$G$, by Corollary~\ref{1.3} and Proposition~\ref{1.4}.

On the other hand, $\sig=p^n$ and $\sig-1=q$ for some primes $p$
and $q$. According to Lemma~\ref{3.9}, either $\sig=3$ or $q$ is a
Mersenne prime. Then $G/\Phi(G)$ is isomorphic either to the group
(i) or to (ii), depending on $\sig=3$ or $\sig\neq 3$. We give a
proof for the case when $\sig\neq 3$. The first case is discussed
similarly.

We have $G/C_{\sig+1}\cong
E_{\sig}\rtimes(M_{\sig+1}/C_{\sig+1})$. Since each maximal
subgroup of $M_{\sig+1}/C_{\sig+1}$ is normal of index $q$,
$M_{\sig+1}/C_{\sig+1}$ is a $q$-group. But
$q^2\nmid|\Aut(E_{\sig})|$, so we have
$|M_{\sig+1}/C_{\sig+1}|=q$. We may assume that $C_{\sig+1}\nleq
M_1$, because $G$ has more than one normal maximal subgroup. Hence
$ G/(M_1\cap C_{\sig+1})\cong E_{\sig}\rtimes\Z_q\times\Z_q.$ Now
since $\f{m}(G/(M_1\cap C_{\sig+1}))=\f{m}(G)$, we have
$\Phi(G)=M_1\cap C_{\sig+1}$.

\vspace*{0.2cm} \noindent {\bf\it Step 2.}\\ We next assume that
$G$ has a non-normal maximal subgroup of index $\sig-1$. By
Proposition 1.3, we may assume that $(M_1,\cdots,M_{\sig})$
satisfies the condition (i) of Theorem~\ref{1.1}. We now
distinguish three following cases:

\vspace*{0.2cm}\noindent {\bf Case 1.} $M_j$ is normal in $G$,
where $\sig+1\leq j\leq 2\sig$.

In this case by Corollary~\ref{1.3}, the indices of $M_j$'s
$(\sig+1\leq j\leq 2\sig)$ are mutually coprime or equal to
$\sig-1$, which lead to the groups (iii) and (iv), respectively.
We give a proof for the latter case. We have
$G/C_2\cong\Z_{\sig-1}\rtimes\Z_t$ and
$G/C\cong\Z_{\sig-1}\times\Z_{\sig-1}$, where $C=M_{\sig+1}\cap
M_{\sig+2}$. Since the non-trivial normal subgroups of $G/C$ are
of indices $\sig-1$ and those of $G/C_2$ have indices less than
$\sig-1$, we have $G=CC_2$. Now the proof is completed by the fact
that $G/(C\cap C_2)\cong G/C\times G/C_2$.

\vspace*{0.2cm} \noindent {\bf Case 2.} There exists an integer
$j$ with $\sig+1\leq j\leq 2\sig$ such that $M_j$ is non-normal of
index $\sig-1$ in $G$.

In this case, we may assume that $\{M_{\sig+1},...,M_{2\sig-1}\}$
is the conjugacy class of $M_j$. Therefore, $M_{2\sig}$ is a
normal subgroup of $G$ with index coprime to that of $M_1$. It
follows, from Proposition 1.4(i), that $M_2/C_2$ and
$M_{\sig+1}/C_{\sig+1}$ are cyclic. If $G=C_2C_{\sig+2}$ then, by
Theorem~\ref{3.4'}(ii), $G/\Phi(G)$ is isomorphic to (v), and
$|M_2/C_2|$, $|M_{\sig+1}/C_{\sig+1}|$ are coprime. Hence in this
case $\sig-1$ is not a Fermat prime.

We now suppose that $G\neq C_2C_{\sig+1}$. If $M_{2\sig}$ contains
$C_2$ or $C_{\sig+1}$, then $G/\Phi(G)$ is isomorphic to the group
(vi), by Theorem~\ref{3.5}. Otherwise, $G/\Phi(G)$ will be
isomorphic to the group (vii).

\vspace*{0.2cm} \noindent {\bf Case 3.} There exists an integer
$j$ with $\sig+1\leq j\leq 2\sig$ such that $M_j$ is non-normal of
index $\sig$ in $G$.

In this case, $\sig$ and $\sig-1$ are prime powers. By
Lemma~\ref{3.9}, one of the following statements holds:\\ (a)\quad
$\sig=q+1$, where $q$ is a Mersenne prime;\\ (b)\quad $\sig$ is a
Fermat prime;\\ (c)\quad $\sig=9$.\\ In what follows, we shall
show that the cases (b) and (c) cannot occur. In the case of (b),
$|M_{\sig+1}/C_{\sig+1}|$ is a power of $2$. However, by
Proposition 1.4(iii), $M_2/C_2$ is cyclic of odd order, as
$O_2(M_2/C_2)=1$. So $G$ must have two normal maximal subgroups, a
contradiction.

\vspace*{0.2cm}Now suppose, by way of contradiction, that the case
(c) occurs. Then $G/C_2\cong$ $E_8\rtimes (M_2/C_2)$, which shows
that $M_2/C_2$ has no maximal subgroup of index $9$. Therefore,
$M_2/C_2$ has just one maximal subgroup, proving that $M_2/C_2$ is
cyclic. In view of Proposition~\ref{3.8}, we have
$M_2/C_2\cong\Z_7$, so that $|G:M_1|=7$. On the other hand,
$G/C_{10}\cong E_9 \rtimes(M_{10}/C_{10})$. So $7\nmid
|M_{10}/C_{10}|$, from which we deduce that the group
$H=M_{10}/C_{10}$ has a maximal subgroup $K$ of index $8$. Now
according to Lemma~\ref{1.2}, $\sig(H)=9$, whence $K/\Core_H(K)$
is cyclic, by Proposition 1.4. By using Proposition~\ref{3.8},
$|K/\Core_H(K)|=7$, which is impossible.

Therefore, the case (a) is left to be considered. In this case we
first show that $M_{\sig+1}/C_{\sig+1}$ has a non-normal maximal
subgroup. Assume that this not the case. So
$M_{\sig+1}/C_{\sig+1}$ is cyclic, whence it has just one maximal
subgroup. Now using Proposition~\ref{3.8}, we see that
$|M_{\sig+1}/C_{\sig+1}|=q$. Thus, $G$ has a normal maximal
subgroup of index $q$. However, the only normal maximal subgroup
of $G$ is of index less than $q$, which is impossible. Hence
$M_{\sig+1}/C_{\sig+1}$ has a non-normal maximal subgroup, and so
$C_{\sig+1}\leq C_2$. Now we take $H=M_{\sig+1}/C_{\sig+1}$. Since
$\f{m}(H)=\sig(H)=\sig$ and $H$ has a non-normal subgroup of index
$\sig-1$, $H/\Phi(H)\cong\Z_q\rtimes\Z_r$, for some $r$ with
$\f{m}(\Z_r)=1$. It follows that $H$ has a unique Sylow
$q$-subgroup by the Frattini argument. But $q^2\nmid
|\Aut(E_{\sig})|$, which shows that $H\cong\Z_q\rtimes\Z_t$, for
some $t\in \N$ with $\f{m}(\Z_t)=1$. Now, since
$C_{\sig+1}=\Phi(G)$ and $G/C_{\sig+1}=E_{\sig}\rtimes H$, we
conclude that $G/\Phi(G)$ is isomorphic to the group (viii). \Bo


\begin{thebibliography}{99}
\bibitem{a}\label{Cohn} J. H. E. Cohn, {\it On n-sum groups},
 Math. Scand., {\bf 75} (1994), 44-58.
\bibitem{b}\label{DoHa} K. Doerk and T. Hawkes, {\it Finite soluble
groups}, ed., Gruyter, Berlin, (1992).
\bibitem{c}\label{JaMo} A. Jamali and H. Mousavi, {\it A note on the $\sig$-covers of finite
soluble groups}, Algebra Colloquium, (to appear).
\bibitem{d}\label{Thes} H. Mousavi,{\it Characterization of some
finite groups with specific subgroups}, Ph.D. thesise, University
for Teacher Education, (1999), 39-64
\bibitem{e}\label{Paz} \"{G}.
Pazderski, {\it Uber maximale untergruppen endlicher gruppen},
Math. Nachr. {\bf 26}, (1963), 307-319.
\bibitem{f}\label{Wang} Wang Jing, {\it The number of maximal
subgroups and their types}, (Chinese. English Summery), Pure Appl.
Math. {\bf 5} (1989), 24-33.
\end{thebibliography}
\end{document}